\documentclass[a4paper]{amsart}
\begin{document}

\newcommand{\ea}{\mbox{{\bf a}}}
\newcommand{\eu}{\mbox{{\bf u}}}
\newcommand{\ueu}{\underline{\eu}}
\newcommand{\ueo}{\overline{u}}
\newcommand{\oeu}{\overline{\eu}}
\newcommand{\ew}{\mbox{{\bf w}}}
\newcommand{\ef}{\mbox{{\bf f}}}
\newcommand{\eF}{\mbox{{\bf F}}}
\newcommand{\eC}{\mbox{{\bf C}}}
\newcommand{\en}{\mbox{{\bf n}}}
\newcommand{\eT}{\mbox{{\bf T}}}
\newcommand{\eL}{\mbox{{\bf L}}}
\newcommand{\eR}{\mbox{{\bf R}}}
\newcommand{\eV}{\mbox{{\bf V}}}
\newcommand{\eU}{\mbox{{\bf U}}}
\newcommand{\ev}{\mbox{{\bf v}}}
\newcommand{\eve}{\mbox{{\bf e}}}
\newcommand{\uev}{\underline{\ev}}
\newcommand{\eY}{\mbox{{\bf Y}}}
\newcommand{\eK}{\mbox{{\bf K}}}
\newcommand{\eP}{\mbox{{\bf P}}}
\newcommand{\eS}{\mbox{{\bf S}}}
\newcommand{\eJ}{\mbox{{\bf J}}}
\newcommand{\eB}{\mbox{{\bf B}}}
\newcommand{\eH}{\mbox{{\bf H}}}
\newcommand{\leb}{\mathcal{ L}^{n}}
\newcommand{\eI}{\mathcal{ I}}
\newcommand{\eE}{\mathcal{ E}}
\newcommand{\hen}{\mathcal{H}^{n-1}}
\newcommand{\eBV}{\mbox{{\bf BV}}}
\newcommand{\eA}{\mbox{{\bf A}}}
\newcommand{\eSBV}{\mbox{{\bf SBV}}}
\newcommand{\eBD}{\mbox{{\bf BD}}}
\newcommand{\eSBD}{\mbox{{\bf SBD}}}
\newcommand{\ecs}{\mbox{{\bf X}}}
\newcommand{\eg}{\mbox{{\bf g}}}
\newcommand{\paromega}{\partial \Omega}
\newcommand{\gau}{\Gamma_{u}}
\newcommand{\gaf}{\Gamma_{f}}
\newcommand{\sig}{{\bf \sigma}}
\newcommand{\gac}{\Gamma_{\mbox{{\bf c}}}}
\newcommand{\deu}{\dot{\eu}}
\newcommand{\dueu}{\underline{\deu}}
\newcommand{\dev}{\dot{\ev}}
\newcommand{\duev}{\underline{\dev}}
\newcommand{\weak}{\rightharpoonup}
\newcommand{\weakdown}{\rightharpoondown}
\renewcommand{\contentsname}{ }

\newtheorem{rema}{Remark}[section]
\newtheorem{thm}{Theorem}[section]
\newtheorem{lema}{Lemma}[section]
\newtheorem{prop}{Proposition}[section]
\newtheorem{cor}{Corollary}[section]
\newtheorem{defi}{Definition}[section]
\newtheorem{conje}{Conjecture}[]
\newtheorem{exempl}{Example}[section]
\newtheorem{opp}{Open Problem}[]
\renewcommand{\contentsname}{ }
\newenvironment{rk}{\begin{rema}  \em}{\end{rema}}
\newenvironment{exemplu}{\begin{exempl}  \em}{\hfill $\surd$
\end{exempl}}

\title{Symplectic, Hofer and sub-Riemannian geometry}
%\date{21.11.2001}

\begin{abstract}
In this note are proposed some connections between symplectic geometry and 
sub-Riemannian geometry on the Heisenberg group. Almost everything can be extended 
to general Carnot-Carath\'eodory groups.
\end{abstract}


\maketitle

\section{Introduction}

The purpose of this note is to make some connection between the 
sub-Riemannian geometry on Carnot-Carath\'eodory groups and symplectic geometry. 

We shall concentrate here on the Heisenberg group, although it is transparent that almost everything can be done on a general Carnot-Carath\'eodory group. Such  generalisations will be the subject of a forthcoming note.

\section{Carnot-Carath\'eodory groups}

A Carnot-Carath\'eodory (CC) group $N$ is a simply connected nilpotent group 
endowed with a one parameter family of dilations 
$\left\{ \delta_{\varepsilon} \mbox{ : } \varepsilon \in (0,+\infty) 
\right\}$. 

We can identify the group with its algebra and so we get a real vector space 
which is a nilpotent Lie algebra and a Lie group. The 
Baker-Campbell-Hausdorff formula stops after a finite number of steps and it 
makes the connection between the Lie bracket and the group operation. 
From now  on we shall implicitly assume the identification between the 
group and the algebra. This is not stranger than the identification 
of the group $(R,+)$ with it's algebra. The name that $N$ should bear 
is not clear; most of the times it will be called "group", just for convenience. 
For the group operation we shall use the multiplicative notation. 

The dilations $\delta_{\varepsilon}$ are group morphisms and algebra morphisms. 
The group $N$ admits the direct sum decomposition: 
$$N \ = \ \sum_{i=1}^{m} V_{i} \ \ , \ \ [V_{1}, V_{i}] = V_{i+1}  \ , \ 
[V_{1}, V_{m}] = 0$$ 
which is uniquely determined by saying that dilations act like this: 
$$x = \sum_{i=1}^{m} x_{i} \ \mapsto \ \delta_{\varepsilon} x \ = \ 
\sum_{i=1}^{m} \varepsilon^{i} x_{i}$$
The number $m$ is called the step of the group and the number 
$$Q \ = \ \sum_{i=1}^{m} i \ dim \ V_{i}$$
is called the homogeneous dimension of the group.

Any Euclidean norm on $V_{1}$ extends to a (left invariant) distance 
on $N$. Indeed, let $\mid \cdot \mid$ be the norm on $V_{1}$ which we extend 
with $+\infty$ outside $V_{1}$. Then the distance between two points 
$x,y \in N$ is 
$$d(x,y) \ = \ \inf\left\{ \int_{0}^{1} \mid c(t)^{-1}\dot{c}(t)\mid 
\mbox{ d}t \mbox{ : } c(0) = x , \ c(1) = y \right\}$$
This is a good definition of a distance function. It is called the Carnot-Caratheodory
distance associated to the left invariant distribution generated by $V_{1}$. Indeed,
consider the subbundle of the tangent bundle defined in the point $x \in N$ by the plane 
$xV_{1}$. Transport, by using the same left translations, 
the Euclidean norm $\mid \cdot \mid$ all over the distribution. A horizontal curve is a
curve which is tangent almost everywhere to the distribution. The length of a
horizontal curve is measured using the Euclidean norm. Because $V_{1}$ Lie
generates $N$, it follows that any two points can be joined by a horizontal curve with
finite length. The Carnot-Caratheodory distance between two points is then the infimum of
the lengths of horizontal curves joining these two points. A manifold endowed with a
distribution which allows the construction of a CC distance is called a sub-Riemannian 
manifold. A CC group $N$, with the left invariant distribution generated by $V_{1}$ is 
an example of a sub-Riemannian manifold.  For an excellent introduction (and more than
this) into the realm of sub-Riemannian geometry, consult Bella\:{\i}che \cite{bell} and 
Gromov \cite{gromo}.

This distance is left 
invariant, it behaves well with respect to dilations,  hence it is 
generated by a "homogeneous norm".

A continuous function  $x \mapsto \mid x \mid$ from $N$ to 
$[0,+\infty)$ is a homogeneous norm if
\begin{enumerate}
\item[(a)] the set $\left\{ x \in N \ : \ \mid x \mid = 1 \right\}$ 
does not contain $0$.
\item[(b)] $\mid x^{-1} \mid = \mid x \mid$ for any $x \in N$. 
\item[(c)] $\mid \delta_{\varepsilon} x \mid = \varepsilon 
\mid x \mid$ for 
any $x \in N$ and $\varepsilon > 0$. 
\end{enumerate}

Any two homogeneous norms are equivalent. Let $\mid \cdot \mid$ be a
homogeneous norm. Then the set 
$\left\{ x \mbox{ : } \mid x \mid = 1 \right\}$ is compact. 
Moreover, there is a constant $C> 0$ such that for any $x,y \in N$ we have: 
$$ \mid xy \mid \ \leq \ C \ ( \mid x \mid + \mid y \mid ) $$
For the particular homogeneous norm 
$$\mid x \mid_{d} \ = \ d(0,x)$$ 
the constant $C$ equals 1. 

The homogeneous norm induced by the distance might be impossible to compute, 
but computable homogeneous norms are aplenty, for example this one: 
$$\mid \sum_{i=1}^{m} x_{i} \mid \ = \ \sum_{i=1}^{m} \mid x_{i} \mid^{1/i}$$

The Hausdorff measure $\mathcal{H}^{Q}$ is proportional to the 
Lebesgue measure on the vector space $N$. Therefore the metric dimension of $N$ is strictly larger than the topological dimension. 

We can be more precise by saying that the Hausdorff measure $\mathcal{H}^{Q}$ is a bi-invariant Haar measure on the group $N$; moreover, this measure is compatible with the dilations: for any $\varepsilon > 0$ and any set $A$ we have: 
$$\mathcal{H}^{Q} (\delta_{\varepsilon} A) \ = \ \varepsilon^{Q} \mathcal{H}^{Q} (A)$$

We give further examples of CC groups:

{\bf (1.)} $R^{n}$ with addition is the only commutative CC group. 

{\bf (2.)} The Heisenberg group is the first non-trivial example. 
This is the group $H(n) \ = \ R^{2n} \times R$ with the operation: 
$$(x,\bar{x}) (y,\bar{y}) \ = \ (x+y, \bar{x} + \bar{y} + 
\frac{1}{2} \omega(x,y) )$$
where $\omega$ is the standard symplectic form on $R^{2n}$. 
The Lie bracket is 
$$[(x,\bar{x}), (y,\bar{y})] \ = \ (0,\omega(x,y))$$
The direct sum decomposition of (the algebra of the) group is: 
$$H(n) \ = \ V + Z \ , \ \ V \ = \ R^{2n} \times \left\{ 0 \right\}
\ , \ \ Z \ = \  \left\{ 0 \right\} \times R$$ 
$Z$ is the center of the algebra, the group has step 2 and
homogeneous dimension $2n+2$. 

{\bf (3.)} H-type groups. These are two step nilpotent Lie groups $N$ 
endowed with an inner product $(\cdot , \cdot)$, such that the
following {\it orthogonal} direct sum decomposition occurs: 
$$N \ = \ V + Z$$
$Z$ is the center of the Lie algebra. Define now the function 
$$J : Z \rightarrow End(V) \ , \ \ (J_{z} x, x') \ = \ (z, [x,x'])$$ 
The group $N$ is of H-type if for any $z \in Z$ we have 
$$J_{z} \circ J_{z} \ = \  - \mid z \mid^{2} \  I$$
From the Baker-Campbell-Hausdorff formula we see that the group
operation is 
$$(x,z) (x', z') \ = \ (x + x', z + z' + \frac{1}{2} [x,x'])$$
These groups appear naturally as the nilpotent part in the Iwasawa
decomposition of a semisimple real group of rank one. (see \cite{fostein})

{\bf (4.)} The last example is the group of $n \times n$ 
upper triangular matrices, which is
nilpotent of step $n-1$.

\section{Pansu differentiability}

A CC group has it's own concept of differentiability.

In Euclidean spaces, given $f: R^{n} \rightarrow R^{m}$ and 
a fixed point $x \in R^{n}$, one considers the difference function: 
$$X \in B(0,1) \subset R^{n} \  \mapsto \ \frac{f(x+ tX) - f(x)}{t} \in R^{m}$$
The convergence of the difference function as $t \rightarrow 0$ in the uniform 
convergence gives raise to the concept of differentiability in it's classical sense. The same convergence, but in measure, leads to approximate differentiability. 
Other topologies might be considered (see Vodop'yanov \cite{vodopis}).

The point is that the difference function can be written using only dilations 
and the group operation. Indeed, for any function between Carnot groups 
$f: G \rightarrow P$,  for  any fixed point $x \in G$ and $t>0$ we  have the difference function: 
$$X \in B(1) \subset G \  \mapsto \ \delta_{t}^{-1} \left(f(x)^{-1}f\left( 
x \delta_{t}X\right) \right) \in P$$
Pansu's differentiability is obtained from uniform convergence of the difference 
function when $t \rightarrow 0$. 

We introduce first the linear functions. 

\begin{defi}
$F:N \rightarrow N$ is linear if 
\begin{enumerate}
\item[(a)] $F$ is a {\it group} morphism, 
\item[(b)] for any $\varepsilon > 0$ $F \circ \delta_{\varepsilon} \
= \ \delta_{\varepsilon} \circ F$. 
\end{enumerate}
\label{dlin}
\end{defi}

The condition (b) means that $F$, seen as an algebra morphism, 
preserves the grading.

\begin{defi}
Let $f: N \rightarrow N$ and $x \in N$. We say that $f$ is 
(Pansu) differentiable in the point $x$ if there is a linear 
function $Df(x): N \rightarrow N$ such that 
$$\sup \left\{ d(F_{t}(y), Df(x)y) \ \mbox{ : } \ y \in B(0,1)
\right\}$$
converges to $0$ when $t \rightarrow 0$. The functions $F_{t}$
are defined by 
$$F_{t} (y) \ = \ \delta_{t}^{-1} \left( f(x)^{-1} f(x
\delta_{t}y)\right)$$
\label{ddif}
\end{defi}

The definition says that $f$ is differentiable at $x$ if the 
sequence of finite differences $F_{t}$ uniformly converges to a
linear map when $t$ tends to $0$.

Theorem 2, Pansu \cite{pansu} contains the Rademacher theorem for Carnot groups: 

\begin{thm}
Let $f: M \rightarrow N$ be a Lipschitz  function between CC groups. 
Then $f$ is differentiable almost everywhere.
\label{ppansu}
\end{thm}

We shall need further a form  of the Poincar\'e inequality for CC groups. 
(This theorem is an adaptation to our needs of a theorem true on metric measure  spaces). 

\begin{thm}
If $\phi$ is Lipschitz then 
$$\inf_{z \in N} \int_{B(x,r)} d(\phi(y), z) \mbox{ d}y \ \leq \ C \int_{B(x,r)}
Lip(D\phi(y)) \mbox{ d}y$$
\label{tpoin}
\end{thm}

Finally, a result showing that in a CC group, Lipschitz functions have Lusin property. 

\begin{thm}
If $\phi: M \rightarrow N$ is Lipschitz then it transports negligible sets into negligible ones. 
\end{thm}

\section{Example: the Heisenberg group}

The Heisenberg group $H(n) = R^{2n+1}$ is a 2-step  nilpotent group with the 
operation: 
$$(x,\bar{x})  (y,\bar{y}) = (x + y, \bar{x} + \bar{y} + \frac{1}{2} \omega(x,y))$$
where $\omega$ is the standard symplectic form on $R^{2n}$. We shall identify 
the Lie algebra with the Lie group. The bracket is 
$$[(x,\bar{x}),(y,\bar{y})] = (0, \omega(x,y))$$
The Heisenberg algebra is generated by 
$$V = R^{2n} \times \left\{ 0 \right\}$$ 
and we have the relation $V + [V,V] = H(n)$.

The group $H(n)$ has dilations 
$$\delta_{\varepsilon} (x,\bar{x}) = (\varepsilon x , \varepsilon^{2} \bar{x})$$

We shall denote by $GL(H(n))$ 
the group of invertible linear transformations and by $SL(H(n))$ the 
subgroup of volume preserving ones. 

\begin{prop}
We have the isomorphisms $$GL(H(n)) = CSp(n) \ , \ \ SL(H(n)) = Sp(n)$$
\label{p1}
\end{prop}

\begin{proof}
By direct computation. We are looking first for the algebra isomorphisms of 
$H(n)$. Let the matrix 
$$\left( \begin{array}{cc} 
          A & b \\
          c & a 
         \end{array} \right)$$
represent such a morphism, with $A \in gl(2n, R)$, $b,c \in R^{2n}$ and 
$a \in R$. The bracket preserving condition reads: for any $(x,\bar{x}), 
(y,\bar{y}) \in H(n)$ we have
$$(0,\omega(A x + \bar{x} b, A y + \bar{y} b)) = ( \omega(x,y) b, a \omega(x,y))$$
We find therefore $b = 0$ and $\omega(Ax, Ay) = a \omega(x,y)$, so $A \in 
CSp(n)$ and $a \geq 0$, $a^{n} = \det A$. 

The preservation of the grading gives $c=0$. The volume preserving condition 
means $a^{n+1} = 1$ hence $a= 1$ and $A \in  Sp(n)$. 
\end{proof}

\subsection{Get acquainted with Pansu differential}
{\bf Derivative of a curve:} 
Let us see which are the smooth (i.e. derivable) curves. Consider 
$\tilde{c}: [0,1] \rightarrow H(n)$, $t \in (0,1)$ and 
$\varepsilon > 0$ sufficiently small. Then the finite difference 
function associated to $c,t, \varepsilon$ is 
$$C_{\varepsilon}(t)(z) \ = \
\delta_{\varepsilon}^{-1} \left( \varepsilon\tilde{c}(t)^{-1}
\tilde{c}(t + \varepsilon z) \right)$$ 
After a short computation we obtain: 
$$C_{\varepsilon}(t)(z) \ = \ \left(\frac{c(t + \varepsilon z) -
c(t)}{\varepsilon} , \frac{\bar{c}(t+\varepsilon z) -
\bar{c}(t)}{\varepsilon^{2}} - \frac{1}{2} \omega(c(t), 
\frac{c(t + \varepsilon z) -
c(t)}{\varepsilon^{2}}) \right)$$
When $\varepsilon \rightarrow 0$ we see that the finite difference
function converges if: 
$$\dot{\bar{c}}(t) \ = \ \frac{1}{2} \omega(c(t), \dot{c}(t))$$
and the limit equals
$$D c(t) z \ = \ z (  \dot{c}(t) , \frac{1}{2} \ddot{\bar{c}}(t) - 
\frac{1}{4} \omega(c(t), \ddot{c}(t))$$
Hence the curve has to be horizontal; in this case we see that 
$$D c(t) z \ = \ z ( \dot{c}(t), 0)$$
This is almost the tangent to the curve. The tangent is obtained 
by taking $z = 1$ and the left translation of $Dc(t) 1$ by $c(t)$.

The horizontality condition implies that,  given a curve 
$t \mapsto c(t) \in R^{2n}$, there is
only one horizontal curve $t \mapsto (c(t),\bar{c}(t))$, such that 
$\bar{c}(0) = 0$. This curve is called the lift of $c$. 

{\bf Derivative of a functional:} 
Take now $f: H(n) \rightarrow R$ and compute its Pansu derivative. 
The finite difference function is 
$$F_{\varepsilon}(x \bar{x})(y,\bar{y}) \ = \ 
\left( f(x + \varepsilon y, \bar{x} + \frac{\varepsilon}{2} \omega(x,y) +
\varepsilon^{2} \bar{y}) - f
(x,\bar{x}) \right)/\varepsilon$$
Suppose that $f$  is (classically) derivable. Then it is Pansu derivable  and this
derivative has the expression: 
$$Df(x , \bar{x}) (y,\bar{y}) \ = \ \frac{\partial f}{\partial
x}(x,\bar{x}) y
+ \frac{1}{2} \omega(x,y) \ \frac{\partial f}{\partial
\bar{x}}(x,\bar{x}) $$
Remark that not any Pansu derivable functional is derivable in the classical sense. As an
example, check that the square of any homogeneous norm is Pansu derivable everywhere, 
but not derivable everywhere in the classical sense.

\section{Symplectomorphisms, capacities and Hofer distance}

Symplectic capacities are invariants of under the action of the symplectomorphisms group. 
Hofer geometry is the geometry of the group of Hamiltonian diffeomorphisms, with respect 
to the Hofer distance. For an introduction into the subject see Hofer, Zehnder
\cite{hozen} chapters 2,3 and 5, and Polterovich \cite{polte}, chapters 1,2. Here is 
a quick introduction of the notions that we shall need further. 

A symplectomorphism is an invertible map $\phi: A \rightarrow B$,  between 
two symplectic manifolds, such that $\phi$ and $\phi^{-1}$ preserve the symplectic 
structure. A symplectic capacity is a map which associates to any symplectic 
manifold $(M,\omega)$ a number $c(M,\omega) \in [0,+\infty]$. Symplectic capacities 
are special cases of conformal symplectic invariants, described by: 
\begin{enumerate}
\item[A1.] Monotonicity: $c(M,\omega) \leq c(N,\tau)$ if there is a symplectic 
embedding  from $M$ to $N$, 
\item[A2.] Conformality: $c(M,\varepsilon \omega) = \mid \varepsilon \mid c(M,\omega)$ 
for any $\alpha \in R$, $\alpha \not = 0$. 
\end{enumerate}

We can see a conformal symplectic invariant from another point of view. Take 
a symplectic manifold $(M,\omega)$ and consider the invariant defined over the class 
of Borel sets $B(M)$,  (seen as embedded submanifolds). In the particular case of 
$R^{2n}$ with the standard symplectic form, an invariant is a function 
$c: B(R^{2n}) \rightarrow [0,+\infty]$ such that: 
\begin{enumerate}
\item[B1.] Monotonicity: $c(M) \leq c(N)$ if there is a symplectomorphism $\phi$ such 
that $\phi(M) \subset N$,
\item[B2.] Conformality: $c(\varepsilon M ) = \varepsilon^{2} c(M)$ for any 
$\varepsilon \in R$. 
\end{enumerate} 

An invariant is nontrivial if it takes finite values on sets with infinite volume, 
like cylinders: 
$$Z(R) = \left\{ x \in R^{2n} \mbox{ : } x^{2}_{1} + x_{2}^{2} < R \right\}$$

There exist highly nontrivial invariants, as the following theorem shows: 

\begin{thm} (Gromov's squeezing theorem) There is a symplectomorphism 
$\phi$ such that $\phi(B(r)) \subset Z(R)$ if and only if $r \leq R$. 
\end{thm}

This theorem permits to define the invariant: 
$$c(A) \ = \ \sup \left\{ R^{2} \mbox{ : } \exists \phi(B(R)) \subset A \right\}$$
called Gromov's capacity. 

Another important invariant is Hofer-Zehnder capacity. In order to introduce this 
we need the notion of a Hamiltonian flow.  

A flow of symplectomorphisms $t \mapsto \phi_{t}$ is Hamiltonian if there is 
a function $H: M \times R \rightarrow R$ such that for any time $t$ and place $x$ we have 
$$\omega( \dot{\phi}_{t}(x), v) \ = \ dH(\phi_{t}(x),t) v$$
for any $v \in T_{\phi_{t}(x)}M$. 

Let $H(R^{2n})$ be the set of compactly supported Hamiltonians. Given 
a set $A \subset R^{2n}$, the class of admissible Hamiltonians is $H(A)$, made by all  compactly supported maps in $A$ such that the generated  Hamiltonian flow does not have closed orbits of periods smaller than 1. Then the Hofer-Zehnder capacity 
is defined by: 
$$hz(A) \ = \ \sup \left\{ \| H \|_{\infty} \mbox{ : } H \in H(A) \right\}$$

Let us denote by $Ham(A)$ the class of Hamiltonian diffeomorphisms compactly supported in $A$. A Hamiltonian diffeomorphism is the time one value of a Hamiltonian flow. In the case which interest us, that is $R^{2n}$, $Ham(A)$ is the connected component of the identity in the group of compactly supported symplectomorphisms. 

A curve of Hamiltonian diffeomorphisms (with compact support) is a Hamiltonian flow. 
For any such curve $t \mapsto c(t)$ we shall denote by $t\mapsto H_{c}(t, \cdot)$ the associated Hamiltonian  function (with compact support). 

On the group of Hamiltonian diffeomorphisms there is a bi-invariant distance introduced by Hofer. This is given by the expression: 
$$d(\phi,\psi) \ = \ \inf \left\{ \int_{0}^{1} \| H_{c}(t) \|_{\infty, R^{2n}} \mbox{ d}t \mbox{ : } c: [0,1] \rightarrow Ham(R^{2n}) \right\}$$

It is easy to check that $d$ is indeed bi-invariant and it satisfies the triangle property. It is a deep result that $d$ is non-degenerate, that is 
$d(id, \phi) \ = \ 0 $ implies $\phi = \ id$. 

With the help of the Hofer distance one can define another symplectic invariant, called displacement energy. For a set $A \subset R^{2n}$ the displacement energy is: 
$$de(A) \ = \ \inf \left\{ d(id, \phi) \ \mbox{ : } \phi \in Ham(R^{2n}) \ \ , \ \phi(A) \cap A = \emptyset \right\}$$

\section{Volume preserving and symplectic diffeomorphisms}

In this section we shall do some computations connected to the 
group of volume preserving diffeomorphisms of $H(n)$. This
computations require some regularity; we shall see in the next
section that we can work without it.

\subsection{Volume preserving diffeomorphisms} 
Let $\tilde{f} = (f,\bar{f}) : H(n) \rightarrow H(n)$ be a smooth 
function (in the sense that $\tilde{f}$ is at least of class $C^{2}$
with respect to the classical differential structure). 
We shall compute: 
$$D \tilde{f} ((x,\bar{x})) (y,\bar{y}) \ = \ 
\lim_{\varepsilon \rightarrow 0} \delta_{\varepsilon^{-1}}  \left( 
\left(\tilde{f}(x,\bar{x})\right)^{-1}  \tilde{f}\left((x,\bar{x})  \delta_{\varepsilon} (y,\bar{y})\right)\right) $$
We know that $D \tilde{f}(x,\bar{x})$ has to be a linear mapping. 

After a short computation we see that we have to pass to the limit 
$\varepsilon \rightarrow 0$ in the following expressions (representing the two 
components of $D \tilde{f} ((x,\bar{x})) (y,\bar{y})$): 
\begin{equation}
\frac{1}{\varepsilon} \left( f\left(x+ \varepsilon y, \bar{x} + \varepsilon^{2} 
\bar{y} + \frac{\varepsilon}{2} \omega(x,y)\right) - f(x,\bar{x}) \right)
\label{exp1}
\end{equation}
\begin{equation}
\frac{1}{\varepsilon^{2}} \left( \bar{f}\left(x+ \varepsilon y, \bar{x} + \varepsilon^{2} 
\bar{y} + \frac{\varepsilon}{2} \omega(x,y)\right) - \bar{f}(x,\bar{x}) - 
\frac{1}{2}\omega\left(f(x,\bar{x}), f\left(x+ \varepsilon y, \bar{x} + \varepsilon^{2} 
\bar{y} + \frac{\varepsilon}{2} \omega(x,y)\right)\right)\right)
\label{exp2}
\end{equation}

For the first component \eqref{exp1} tends to 
$$\frac{\partial f}{\partial x} (x,\bar{x}) y + \frac{1}{2} \frac{\partial f}{\partial \bar{x}} (x,\bar{x}) \omega(x,y)$$
The terms of order $\varepsilon$ must cancel in the second component \eqref{exp2}. We obtain the cancelation 
condition (we shall omit from now on the argument $(x,\bar{x})$ of all functions): 
\begin{equation}
\frac{1}{2} \omega(x,y) \frac{\partial \bar{f}}{\partial \bar{x}} - 
\frac{1}{2} \omega(f, \frac{\partial f}{\partial x} y) - 
\frac{1}{4} \omega(x,y) \omega(f, \frac{\partial f}{\partial \bar{x}}) + 
\frac{\partial \bar{f}}{\partial x} \cdot y \ = \ 0
\label{cancel}
\end{equation}
The second component tends to 
$$\frac{\partial \bar{f}}{\partial \bar{x}} \bar{y} - \frac{1}{2} \omega(f, 
\frac{\partial f}{\partial \bar{x}}) \bar{y}$$
The tangent Lie algebra morphism to the group morphism $D \tilde{f}(x,\bar{x})$ 
is the matrix: 
\begin{equation}
d \tilde{f}(x,\bar{x}) \ = \ \left( \begin{array}{cc}
\frac{\partial f}{\partial x} + \frac{1}{2} \frac{\partial f}{\partial \bar{x}} 
\otimes Jx & 0 \\
0 & \frac{\partial \bar{f}}{\partial \bar{x}} - \frac{1}{2} \omega(f, 
\frac{\partial f}{\partial \bar{x}}) 
       \end{array} \right)
\label{tang}
\end{equation}
We shall suppose now that $\tilde{f}$ is volume preserving. According to 
proposition \ref{p1}, this means: 
\begin{equation}
\frac{\partial f}{\partial x} + \frac{1}{2} \frac{\partial f}{\partial \bar{x}} 
\otimes Jx \ \in Sp(n) 
\label{c1}
\end{equation}
\begin{equation}
\frac{\partial \bar{f}}{\partial \bar{x}} - \frac{1}{2} \omega(f, 
\frac{\partial f}{\partial \bar{x}}) = 1 
\label{c2}
\end{equation}
The cancelation condition \eqref{cancel} and \eqref{c2} give
\begin{equation}
\frac{\partial \bar{f}}{\partial x} y \ = \  \frac{1}{2} \omega(f,\frac{\partial 
f}{\partial x} y ) \ - \ \frac{1}{2} \omega(x,y)
\label{c3}
\end{equation}

These conditions describe completely the class of volume preserving diffeomorphisms 
of $H(n)$. Conditions \eqref{c2} and \eqref{c3} are in fact differential equations 
for the function $\bar{f}$ when $f$ is given. However, there is a compatibility 
condition in terms of $f$ which has to be fulfilled for  \eqref{c3} to have 
a solution $\bar{f}$. Let us look closer to \eqref{c3}. We can see the symplectic 
form $\omega$ as a closed 2-form. Let $\lambda$ be a 1-form such that 
$d \lambda = \omega$. If we take the (regular) differential with respect 
to $x$ in \eqref{c3} we quickly obtain the compatibility condition
\begin{equation}
\frac{\partial f}{\partial x} \ \in \ Sp(n)
\label{c4}
\end{equation}
and \eqref{c3} takes the form: 
\begin{equation}
2 \ d \bar{f} \ = \ f^{*} \lambda \ - \ \lambda
\label{c5}
\end{equation}
(all functions seen as functions of $x$ only).

Conditions \eqref{c4} and \eqref{c1} imply: there is a scalar function 
$\mu = \mu(x,\bar{x})$ such that 
$$\frac{\partial f}{\partial \bar{x}} \ = \ \mu \ Jx $$
Let us see what we have until now: 
\begin{equation}
\frac{\partial f}{\partial x} \ \in \ Sp(n) 
\label{cc1}
\end{equation}
\begin{equation}
\frac{\partial \bar{f}}{\partial x} \ = \ \frac{1}{2} \left[ \left( 
\frac{\partial f}{\partial x}\right)^{T} J f \ - \ J x \right]
\label{cc2}
\end{equation}
\begin{equation}
\frac{\partial \bar{f}}{\partial \bar{x}} \ = \ 1 + \frac{1}{2} \omega(f, 
\frac{\partial f}{\partial \bar{x}}) 
\label{cc3}
\end{equation}
\begin{equation}
\frac{\partial f}{\partial \bar{x}} \ = \ \mu \ Jx
\label{cc4}
\end{equation}
Now, differentiate \eqref{cc2} with respect to $\bar{x}$ and use \eqref{cc4}. In the same time 
differentiate \eqref{cc3} with respect to $x$. From the equality 
$$\frac{\partial^{2} \bar{f}}{\partial x \partial \bar{x}} \ = \ 
\frac{\partial^{2} \bar{f}}{\partial \bar{x} \partial x}$$ 
we shall obtain by straightforward computation $\mu = 0$. 

In order to properly formulate the conclusion of our computations, we need 
the following definition: 

\begin{defi}
$Diff^{2}(H(n),vol)$ is  the group of volume preserving 
diffeomorphisms $\tilde{\phi}$ of $H(n)$ such that 
$\tilde{\phi}$ and it's inverse have (classical) regularity $C^{2}$. 
In the same way we define $Sympl^{2}(R^{2n})$ to be the group of 
$C^{2}$ symplectomorphisms of $R^{2n}$. 
\end{defi}

We have  proven the following result: 

\begin{prop}
Let $Diff^{2}(H(n),vol)$ be the group of volume preserving 
diffeomorphisms 
of $H(n)$  and $Sympl^{2}(R^{2n})$ be the group of symplectomorphisms of 
$R^{2n}$, both with previous mentioned regularity. 
Then we have the isomorphism of groups
$$Diff(H(n),vol) \ = \ Sympl(R^{2n}) \times R$$
given by the mapping 
$$\tilde{f} = (f,\bar{f}) \  \in \ Diff(H(n),vol) \ \mapsto \ \left( 
f \in Sympl(R^{2n}) , \bar{f}(0,0) \right)$$
The inverse of this isomorphism has the expression
$$\left( f \in Sympl(R^{2n}) , a \in R \right) \ \mapsto  \  \tilde{f} = (f,\bar{f}) \  \in \ Diff(H(n),vol)$$ 
$$\tilde{f}(x,\bar{x}) \ = \ (f(x), \ \bar{x} + F(x))$$
where $F(0)= a$ and $dF \ = \ f^{*} \lambda \ - \ \lambda$. 
\label{t1}
\end{prop}

\subsection{Hamilton's equations}

Let $A \subset R^{2n}$ be a set. $Sympl(A)_{c}$ is the group of symplectomorphisms 
with compact support in $A$, that is the group of all symplectomorphisms which differ 
from the identity map on a compact set included in $A$.

\begin{defi}
For any flow $t \mapsto \phi_{t} \in Sympl(A)_{c}$
Denote by $\phi^{h}(\cdot, x))$ the  horizontal flow in 
$H(n)$ obtained by the lift of all curves $t \mapsto \phi(t,x)$ and 
by $\tilde{\phi}(\cdot, t)$ the flow obtained by the lift of all 
$\phi_{t}$. The vertical flow is defined by the expression 
\begin{equation}
\phi^{v} \ = \ \tilde{\phi}^{-1} \circ \phi^{h}
\label{hameq}
\end{equation}
\end{defi}

Relation \eqref{hameq} can be seen as Hamilton equation. 

\begin{prop}
Let $t \in [0,1] \mapsto \phi^{v}_{t}$ be a curve of diffeomorphisms 
of $H(n)$ satisfying the equation: 
\begin{equation}
\frac{d}{dt} \phi^{v}_{t}(x,\bar{x}) \ = \ (0, H(t,x)) \ \ , \ \ 
\phi^{v}_{0} \ = \ id_{H(n)}
\label{ham}
\end{equation}
Then the flow $t \mapsto \phi_{t}$ which satisfies \eqref{hameq} and 
$\phi_{0} \ = \ id_{R^{2n}}$ is the Hamiltonian flow generated by
$H$. 

Conversely, for any Hamiltonian flow $t \mapsto \phi_{t}$, generated 
by $H$, the vertical flow $ t \mapsto \phi^{v}_{t}$ satisfies the
equation \eqref{ham}. 
\label{pham}
\end{prop}

\begin{proof}
Write the lifts $\tilde{\phi}_{t}$ and $\phi^{h}_{t}$, compute then 
the differential of the quantity 
$\dot{\tilde{\phi}}_{t} - \dot{\phi}^{h}_{t}$ and show that it
equals the differential of $H$. 
\end{proof}

\subsection{Flows and gradients}

We want to know if there is any nontrivial smooth (according to Pansu differentiability) 
flow of volume preserving diffeomorphisms. 

\begin{prop}
Suppose that $t \mapsto \phi_{t} \in Diff^{2}(H(n),vol)$ is a flow such that 
\begin{enumerate}
\item[-] is $C^{2}$ in the classical sense with respect to $(x,t)$, 
\item[-] is horizontal, that is $t \mapsto \phi_{t}(x)$ is a horizontal curve 
for any $x$. 
\end{enumerate}
Then the flow is constant. 
\label{pho}
\end{prop}

\begin{proof}
By direct computation, involving second order derivatives. 
\end{proof}

One should expect such a result to be true, based on two remarks. The first: take a 
flow of left translations, that is a flow $t \mapsto \phi_{t}(x) \ = \ x_{t} x$. 
We can see directly that each $\phi_{t}$ is smooth, because the distribution is 
left invariant. But the flow is not horizontal, because the distribution is not 
right invariant.  The second remark: the flow correspond to a Hamiltonian flow 
with null Hamiltonian function, hence the flow is constant. 

At a first glance it is  dissapointing to see that the group of volume preserving 
diffeomorphisms contains no smooth paths according to the intrinsic calculus 
on CC groups. But this makes the richness of such groups of diffeomorphisms, 
as we shall see.

\section{Rigidity of Lipschitz maps}

We know now that Rademacher theorem is true on CC groups, in the form 
\ref{ppansu}.  Let us define then the class $Diff(H(n), vol, Lip)$ the class of 
locally bi-Lipschitz, volume preserving diffeomorphisms of $H(n)$. We shall need 
the group $Sympl(R^{2n}, Lip)$ of locally bi-Lipschitz symplectomorphisms of $R^{2n}$.

\begin{prop}
Let $\tilde{\phi} \in Diff(H(n), vol, Lip)$. Then there exist 
$\phi \in Sympl(R^{2n}, Lip)$ and $F: R^{2n} \rightarrow R$ such that 
for almost any point $(x,\bar{x}) \in H(n)$ we have: 
$$\tilde{\phi}(x,\bar{x}) \ = \ (\phi(x), \bar{x} + F(x))$$
\label{pn1}
\end{prop}

\begin{proof}
By theorem \ref{ppansu} $\tilde{\phi}$ is almost everywhere derivable and 
the derivative can be written in the particular form: 
$$\left( \begin{array}{cc}
A & 0 \\
0 & b
       \end{array} \right)$$ such that $b^{n} = \det A$ and $A \in CSp(n,R)$. The volume 
preserving condition gives $A \in Sp(n,R)$, $b = 1$. Take now the function 
$$f(x,\bar{x}) \ = \ (\phi(x,\bar{x}), \bar{\phi}(x,\bar{x}) - \bar{x})$$
which is Lipschitz. From the Poincare inequality on lines we see that 
$f(x,\bar{x}) = (\phi(x), F(x))$ with $\phi$, $F$ from the conclusion of the theorem. 
\end{proof}

The use of  Rademacher theorem \ref{ppansu} gives us very quick a result which resembles 
a lot with proposition \ref{t1}. 

Consider now the group $N \times R$ with the group operation defined component wise. 
This is also a CC group. Indeed, consider the family of dilations 
$$\delta_{\varepsilon}(x,t) \ = \ (\delta_{\varepsilon}(x), \varepsilon t)$$
which gives to $N \times R$ the structure of a CC group. The left invariant distribution 
on the group which generates the distance is (the left translation of)  $W_{1} = V_{1} 
\times R$.

An easy proposition is: 

\begin{prop}
Let $N$ be a noncomutative CC group which admits the orthogonal decomposition 
$$N \ = \ V_{1} + [N,N]$$ and satisfies the condition 
$$V_{1} \cap Z(N) \ = \ 0$$  The group of linear transformations of $N\times R$ is then 
$$GL(N \times R) \ = \ \left\{ \left( 
\begin{array}{cc}
A & 0 \\
c & d
       \end{array} \right) \ \mbox{ : } A \in GL(N) \ , \ c \in V_{1} \ , \ d \in R 
\right\}$$
\label{pglnr}
\end{prop}

\begin{proof}
We shall proceed as in the proof of proposition \ref{p1}. We are looking first at the Lie algebra isomorphisms of $N \times R$, with general form 
$$\left( 
\begin{array}{cc}
A & b \\
c & d
       \end{array} \right) $$
We obtain the conditions: 
\begin{enumerate}
\item[(i)] $c$ orthogonal on $[N,N]$, 
\item[(ii)] $b$ commutes with the image of $A$:  [b,Ay] = 0, for any $y \in N$, 
\item[(iii)] $A$ is an algebra isomorphism of $N$. 
\end{enumerate}
From (ii), (iii) we deduce that $b$ is in the center  of $N$ and from (i) we see that $c \in V_{1}$. 

We want now  the isomorphism to commute with dilations. This condition gives: 
\begin{enumerate}
\item[(iv)] $b \in V_{1}$, 
\item[(v)] $A$ commutes with the dilations of $N$. 
\end{enumerate}
(iii) and (v) imply that $A \in GL(N)$ and (iv) that $b = 0$. 
\end{proof}

The analog of proposition \ref{pho} is the following: 

\begin{prop}
Let $t \mapsto \phi_{t} \in Diff(H(n), vol, Lip)$ be a curve such that  $(x,t) \mapsto 
\Phi(x,t) = (\phi_{t}(x), t)$ is a Lipschitz map from $N \times R$ to itself. 
Then $t \mapsto \phi_{t}$ is a constant curve. 
\label{tn2}
\end{prop}

\begin{proof}
By Rademacher theorem \ref{ppansu} for the group $N \times R$ we obtain that 
$\Phi$ is almost everywhere derivable. Use now proposition \ref{pglnr} to deduce the claim. 
\end{proof}

\section{Invariants of volume preserving diffeomorphisms}

From Proposition \ref{pn1} we see that any volume preserving diffeomorphisms 
$\tilde{\phi} \in Diff(H(n), vol, Lip)$ preserves more than the volume measure. 
Indeed, for any bounded Borel set $\tilde{A} \subset H(n)$, denote by $A$ the projection 
on $R^{2n}$. The arbitrary diffeomorphism $\tilde{\phi}$ has the 
form: 
$$ \tilde{\phi}(x,\bar{x}) \ = \ (\phi(x), \bar{x} + F(x))$$
The following quantities are then conserved:  
$$vol( \tilde{\phi}(\tilde{A}) ) \ = \ vol( \tilde{A}) $$
$$vol( \phi(A) ) \ = \ vol( A)$$

We shall denote by $\chi(\tilde{A})$ the weight of $\tilde{A}$, defined by: 
$$\chi(\tilde{A}) \ = \ \frac{vol( \tilde{A})}{vol( A)}$$ 
An easy computation shows that 
$$\chi(\delta_{\varepsilon} \tilde{A}) \ = \ \varepsilon^{2} \chi(\tilde{A})$$

Here are other definitions of volume preserving invariants in the Heisenberg group. 

\begin{defi}
\begin{enumerate}
\item[1.] Take any flow of volume preserving diffeomorphisms 
$t \mapsto \tilde{\phi}_{t}$ such that $A \cap \phi_{1}(A) =
\emptyset$. 
The displacement energy of $\tilde{A}$ is the infimum of the number 
$$ \int_{0}^{1} \| \dot{\phi}^{v}_{t}\|_{\infty, A} \mbox{ d}t$$ 
denoted by $de(\tilde{A})$.  
\item[2.] Two sets are equivalent ($\tilde{A} \equiv \tilde{B}$) if there is a
volume preserving diffeomorphism $\tilde{\phi}$ of $H(n)$ such that 
$\tilde{B} \ =\ \tilde{\phi}(\tilde{A})$. Define then:  
$$ \inf \left\{ \frac{(diam (\tilde{B}))^{2n+2}}{vol(A)} \ \mbox{ : } \ 
\tilde{B} \equiv \tilde{A} \right\}$$
We call this
type of invariant an isodiameter and call it $isod(A)$. 
\item[3.] In the same way one can define an isoperimeter. Proceed as
in 2., but replace the diameter with the perimeter (for the proper definition of the perimeter in the Heinsenberg group see \cite{fraseca}): 
$$ isop(\tilde{A})^{2n+1} = \inf \left\{ 
\frac{(Per(\tilde{B}))^{2n+2}}{vol(A)^{2n+1}} \ \mbox{ : } \ 
\tilde{B} \equiv \tilde{A} \right\}$$
\end{enumerate}
\end{defi}

\section{Hamiltonian diffeomorphisms: more structure}

In this section we look closer to the structure of the group of volume preserving diffeomorphisms of the Heisenberg group.

We shall begin by considering the Riemannian structure on the CC group N. This means that we shall take the group of dilations 
$$ x \in N \mapsto \varepsilon x \in N$$
for all $\varepsilon > 0$. The left invariant distribution on the group is generated by the whole algebra, therefore after choosing a scalar product on $N$ we have a Riemannian manifold. The theory of Pansu derivative applies here and we define the group $Diff(N, vol, bracket)$ of volume preserving local bi-Lipschitz homeomorphisms which are smooth with respect to the Pansu derivative. These are, by the Lusin 
theorem, up to a set of arbitrary small measure, volume preserving diffeomorphisms of $N$ as manifold, such that the differential is not only linear, but also bracket preserving. We have therefore the following proposition: 

\begin{prop}
Any map $\tilde{\phi} \in Diff(N, vol, bracket)$ has a.e. the form: 
$$\tilde{\phi}(x,\bar{x}) \ = \  (\phi(x), \bar{x} + F(x))$$
where $\phi$ is a locally bi-Lipschitz diffeomorphism and $F: R^{2n} \rightarrow R$ 
is locally Lipschitz. 
\label{pn2}
\end{prop}

This group contains two privileged subgroups: 
$$Diff(H(n),vol) \ = \  Diff(H(n), vol, bracket) \cap Diff(H(n), vol, Lip)$$ 
$$Diff(H(n),vert) \ =  \ \left\{ \tilde{\phi} \in Diff(H(n), vol, bracket) \mbox{ :  } \tilde{\phi}(x,\bar{x}) \ = \ (x, \bar{x} + F(x) ) \right\}$$

Take a one parameter subgroup $t \mapsto \tilde{\phi}_{t} \in Diff(H(n),vol)$. We know that it cannot be smooth as a curve in $Diff(H(n),vol)$, but we also know that there are vertical and horizontal flows $t \mapsto \phi_{t}^{v} \ , \ \phi_{t}^{h}$ such that 
we have the decomposition $\tilde{\phi}_{t} \circ \phi_{t}^{v} \ =  \ \phi_{t}^{h}$.  
Unfortunately none of the flows $t \mapsto \phi_{t}^{v} \ , \ \phi_{t}^{h}$ are one parameter groups.

There is more structure here that it seems. Consider the 
class 
$$HAM(H(n)) = Diff(H(n), vol) \times Diff(H(n),vert)$$
For any pair in this class we shall use the notation 
$$(\tilde{\phi}, \phi^{v}) \in Diff(H(n), vol) \times 
Diff(H(n),vert) = HAM(H(n))$$ 
This class forms a group with the operation: 
$$(\tilde{\phi}, \phi^{v}) (\tilde{\psi}, \psi^{v}) \ = \ 
( \tilde{\phi} \circ \tilde{\psi}, \phi^{v} \circ \tilde{\phi} \circ 
\psi^{v} \circ \tilde{\phi}^{-1})$$

\begin{prop}
If $t \mapsto \tilde{\phi}_{t} \in Diff(H(n),vol)$ is an one parameter group then 
$t \mapsto (\tilde{\phi_{t}}, \phi_{t}^{v}) \in HAM(H(n))$ is an one parameter group. 
\label{popg}
\end{prop}

\begin{proof}
The check is left to the reader. Use Definition and Proposition 1, chapter 5, Hofer \& 
Zehnder \cite{hozen}, page 144. 
\end{proof}

We can put a distribution on the group HAM(H(n)) such that a
 curve $t \mapsto (\tilde{\phi_{t}}, \phi^{v}_{t})$ which corresponds 
to a Hamiltonian flow to be a horizontal curve in the group
$HAM(n)$. Indeed, one can define for any pair $(\tilde{\phi},
\phi^{v}) \in HAM(H(n))$ the function 
$$\phi^{h} \ = \ \tilde{\phi} \circ \phi^{v}$$
and then introduce the distribution on $HAM(H(n))$ by requiring that 
a curve $$t \mapsto (\tilde{\phi}_{t}, \phi_{t}^{v}) \in HAM(H(n))$$  
is horizontal if 
$$t \mapsto \phi^{h}_{t}(x,\bar{x})$$ is horizontal for any 
$(x,\bar{x}) \in H(n)$. 

Let $t \mapsto (\tilde{\phi}_{t}, \phi^{v}_{t}) \in HAM(H(n))$ be a curve 
in the group $HAM(H(n))$. Using the general form of a transformation of 
$HAM(H(n))$, we can write the curve like this:
$$t \mapsto \left((\phi_{t}(x), \bar{x} + F_{t}(x)) \ , \ (x, \bar{x} + \Lambda_{t}(x))\right)$$
The tangent space at $(\tilde{\phi} , \phi^{v}) = (\tilde{\phi}_{0}, \phi^{v}_{0})$ 
is parametrized by $\dot{\phi}, \dot{F}, \dot{\Lambda}$. The distribution mentioned before, or the condition for $\phi_{t}^{h}$ to be horizontal, becomes: 
$$\dot{\Lambda}(x) + \dot{F}(x) \ = \ \frac{1}{2} \omega(\phi(x), \dot{phi}(x))$$

Moreover, we can introduce the following length function for horizontal curves: 
$$L \left(t \mapsto (\tilde{\phi}_{t}, \phi_{t}^{v})\right) \ = \ 
\int_{0}^{1} \| \dot{\phi_{t}^{v}} \|_{L^{\infty}(R^{2n})} \mbox{ d}t $$
With the help of the length we are able to endow the group $HAM(H(n))$ as a 
path metric space. The distance is defined by: 
$$dist\left((\tilde{\phi}_{A}, \phi^{v}_{A}), (\tilde{\phi}_{B}, \phi^{v}_{B})\right) 
\ = \ \inf  L\left(t \mapsto (\tilde{\phi}_{t}, \phi_{t}^{v})\right)$$
over all horizontal curves  $t \mapsto (\tilde{\phi}_{t}, \phi_{t}^{v})$ such that 
$$(\tilde{\phi}_{0}, \phi_{0}^{v}) \ = \ (\tilde{\phi}_{A}, \phi_{A}^{v})$$
$$(\tilde{\phi}_{1}, \phi_{1}^{v}) \ = \ (\tilde{\phi}_{B}, \phi_{B}^{v})$$
The distance is not defined for any two points in $HAM(H(n))$. In principle 
the distance can be degenerated.

The group $HAM(H(n))$ acts on $H(n)$ by: 
$$(\tilde{\phi}, \phi^{v}) (x, \bar{x}) \ = \ \phi^{v} \circ \tilde{\phi} (x, \bar{x})$$

Remark  that 
$$Diff(H(n), vert) \ \equiv \  \left\{ id \right\} \times Diff(H(n),vert)$$ 
 is normal in $HAM(H(n))$. We denote by $HAM(H(n))/Diff(H(n),vert)$ the factor group. 

Consider now the natural action of $Diff(H(n),vert)$ on $H(n)$. The space of orbits 
will be denoted by $H(n)/Diff(H(n),vert)$. 

Factorize now the action of $HAM(H(n))$ on $H(n)$, by the group 
$Diff(H(n),vert)$. We obtain an action of $HAM(H(n))/Diff(H(n),vert)$ on $H(n)/Diff(H(n),vert)$. Close, but elementary, examination shows this is exact the action of the symplectomorphisms group on $R^{2n}$.
All in all we have the following 

\begin{thm}
The action of $HAM(H(n))_{c}$ on $H(n)$ descends after reduction with $Diff(H(n), vert)_{c}$ 
to the action of symplectomorphisms group with compact support on $R^{2n} = H(n)/ Diff(H(n), vert)$. 
The distance $dist$ on $HAM(H(n))$ descends to the Hoofer distance on the connected 
component of identity of the symplectomorphisms group. Any volume preserving invariant 
descends to a symplectic invariant. 
\label{tmai}
\end{thm}

\end{document}